\documentclass[12pt,a4paper]{article}
\usepackage{mathptmx}
\usepackage{courier}
\usepackage[dvips]{graphicx}
\usepackage{psfrag}
\usepackage{amsmath,amssymb}        
\usepackage{latexsym}
\usepackage{url}
\usepackage[light]{kpfonts}

\def\ds{\displaystyle}
\def\e{\epsilon}

\def\D{\ensuremath{D}}
\def\FBF{free boundary formulation}

\title{Free Boundary Formulations for two Extended Blasius Problems}
\author{Riccardo Fazio\\
Department of Mathematics, Computer Science,\\ Physical Sciences and Earth Sciences,\\
University of Messina \\
Viale F. Stagno D'Alcontres, 31 \\
98166 Messina, Italy \\
E-mails: \url{rfazio@unime.it} \\
Home-page: \url {http://mat521.unime.it/fazio}} 
\date{\today}
\begin{document}
\maketitle
\begin{abstract}
In this paper we have defined the \FBF \ for two extended Blasius problems.
These problems are of interest in boundary layer theory and are deduced from the governing partial differential equations by using appropriate similarity variables.
The computed results, for the so-called missing initial condition, are favourably compared with recent results available in the literature.
\end{abstract}
\bigskip

\bigskip


\section{Introduction}
In this paper, we define a \FBF \ for two extended Blasius problems available in literature.
These two boundary value problems (BVPs), which will be recalled in the next section, are both defined on a semi-infinite interval.
From a numerical viewpoint the application of the asymptotic boundary condition cannot be enforced in a simple way.
To overcome this drawback several approaches have been studied.

The classical approach for such a condition is to replace it with the same condition prescribed at a (finite) truncated boundary, as described by Fox \cite[p.92]{Fox} or Collatz \cite[pp. 150-151]{Collatz}.
In many cases, this simple approach, used by trial and errors, result to be sufficiently accurate, although sometimes it provides good results only for very large values of the truncated boundary.
Seldom a simple estimate of the error due to the introduced truncated boundary is available, and for instance, in the case of the Blasius problem an analysis based on the scaling properties of the mathematical model is developed by Rubel \cite{Rubel:1955:EET}.

A better approach, provided an asymptotic analysis to find the the appropriate boundary conditions to be imposed at the truncated boundary can be developed, was proposed by de Hoog and Weiss \cite{deHoog:1980:ATB}, Lentini and Keller \cite{Lentini:BVP:1980} and Marcowich \cite{Markowich:TAS:1982,Markowich:ABV:1983}.
Since the imposed conditions are related to the asymptotic behaviour of the solution, then, usually, the obtained numerical results are more accurate that those from the previous approach.
That is, smaller values of the truncated boundary are necessary in this approach compared with the values required by the classical approach.

The \FBF \ main idea is simple to explain: we replace the asymptotic boundary conditions with two boundary conditions given at an unknown free boundary that has to be determined as part of the solution.
This idea was formulated for the first time for the numerical solution of the Blasius problem by Fazio \cite{Fazio:1992:BPF}.
Moreover, its application has been proposed for the problems in boundary layer theory, see Fazio \cite{Fazio:1996:NAN}.

\section{Two extended Blasius problems}
The first extended Blasius problem was already considered by Schowalter \cite{Schowalter:1960:ABL}, Lee and Ames \cite{Lee:1966:SSN}, Lin and Chern \cite{Lin:1979:LBL}, Kim et al. \cite{Kim:1983:MHT}, or Akcay and Y\"ukselen \cite{Akcay:1999:DRN}.
\begin{align}\label{eq:ExBlasius1} 
& {\displaystyle \frac{d^3 f}{d \eta^3}} {\displaystyle \frac{d^{2}f}{d\eta^2}}^{(P-1)} + \frac{1}{2} \; f
{\displaystyle \frac{d^{2}f}{d\eta^2}} = 0 \nonumber \\[-1ex]
&\\[-1ex]
& f(0) = {\displaystyle \frac{df}{d\eta}}(0) = 0, \qquad
{\displaystyle \frac{df}{d\eta}}(\eta) \rightarrow 1 \quad \mbox{as}
\quad \eta \rightarrow \infty \ , \nonumber 
\end{align}
where $f$ and $\eta$ are appropriate similarity variables and $P$ verifies the conditions $1 \le P < 2$.
Liao \cite{Liao:2005:CNP} has found analytically that the extended Blasius problem (\ref{eq:ExBlasius1}) for $P = 2$ admit an infinite number of solutions and, therefore, in his opinion can be considered as a challenging problem for numerical techniques.
Let us remark here, that the case $P=1$ is the classical Blasius problem, see Blasius \cite{Blasius:1908:GFK}.

As far as the second extended Blasius problem is concerned, the mathematical model arises in the study of a 2D laminar boundary-layer with power-law viscosity for Newtonian fluids and is given by, see Schlichting and Gersten \cite{Schlichting:2000:BLT} or Benlahsen et al. \cite{Benlahsen:2008:GBE},
\begin{align}\label{eq:ExBlasius2} 
& {\displaystyle \frac{d}{d\eta}\left(\left|\frac{d^{2}f}{d\eta^2}\right|^{P-1} \frac{d^{2}f}{d\eta^2}\right)} + \frac{1}{P+1} f {\displaystyle \frac{d^{2}f}{d\eta^2}} = 0 \nonumber \\[-1ex]
&\\[-1ex]
& f(0) = {\displaystyle \frac{df}{d\eta}}(0) = 0 \ , \qquad {\displaystyle \frac{df}{d\eta}}(\eta) \rightarrow 1 \quad \mbox{as} \quad \eta \rightarrow \infty \ , \nonumber 
\end{align}
where $f(\eta)$ is the non-dimensional stream function and $P$ is a given positive value bigger than zero.
Let us remark here, that when  $P = 1$ the BVP (\ref{eq:ExBlasius2}) reduces to the celebrated Blasius problem, see Blasius \cite{Blasius:1908:GFK}.

\section{The \FBF \ idea}
As mentioned before, in a \FBF \ we replace the asymptotic boundary conditions with two boundary conditions fixed at an unknown free boundary that has to be determined as part of the solution.
For both the extended Blasius problems considered in this paper, a \FBF \ can be defined in the same way.
We replace the asymptotic boundary condition with the two conditions
\begin{equation}
{\displaystyle \frac{df}{d\eta}}(\eta_\e) \ , \qquad  {\displaystyle \frac{d^{2}f}{d\eta^2}} = \e \ ,
\end{equation}
where $\eta_\e$ is an unknown free boundary and $\e$ is an assigned small value.
Of course, we can verify if $\eta_\e$ goes to infinity as $\e$ goes to zero.
In order to have a valid formulation, this should always be true.  

\section{Numerical results}
In this section, we report the computed numerical results for the \FBF s of the two extended Blasius problem.

\subsection{First extended Blasius}
Before considering to solve numerically the \FBF \ of the first extended Blasius problem it may be convenient to reformulate it in normal form, sse Asher and Russel \cite{Ascher:1981:RBV}.
To this end, we introduce the new variables $ \theta = \eta/\eta_\e$, $u_1 = f(\eta)$, $u_2 = \frac{df}{d\eta}(\eta)$,$u_3 = \frac{d^2f}{d\eta^2}(\eta)$ and $u_4 =\eta_\e$.
So that, the \FBF \ for the problem (\ref{eq:ExBlasius1}) is given by:
\begin{align}\label{eq:ExB1FBFNF} 
& {\displaystyle \frac{du_1}{d\theta}} = \eta_\e u_2 \nonumber \\
& {\displaystyle \frac{du_2}{d\theta}} = \eta_\e u_3 \nonumber \\
& {\displaystyle \frac{du_3}{d\theta}} = - \eta_\e \frac{1}{2} u_1 u_3^{2-P}  \nonumber \\[-1.5ex]
& \\[-1.5ex]
& {\displaystyle \frac{du_4}{d\theta}} = 0 \nonumber\\
& f(0) = {\displaystyle \frac{df}{d\theta}}(0) = 0, \qquad
{\displaystyle \frac{df}{d\theta}}(1) = 1 \quad {\displaystyle \frac{d^2f}{d\theta^2}}(1) = \e \ , \nonumber 
\end{align}
For the numerical solution of the BVP (\ref{eq:ExB1FBFNF}) we used the \textit{bvp4c.m} MATLAB routine with initial iterate given by
\begin{equation}\label{eq:IT}
u_1 = \theta \ , \quad u_2 = 2+\theta \ , \quad u_3 = \theta \ , \quad u_4 = 1 \ .
\end{equation}
As it is easily seen, these are coarse approximations of the actual solution components.

\begin{figure}[!hbt]
\centering
\psfrag{e}[1][]{$\eta$}  
\psfrag{f1}[l][]{$\ds \frac{df}{d\eta}$} 
\psfrag{f2}[l][]{$\ds \frac{d^2f}{d\eta^2}$} 
\psfrag{f}[l][]{$f$} 
\includegraphics[width=14cm,height=14cm]{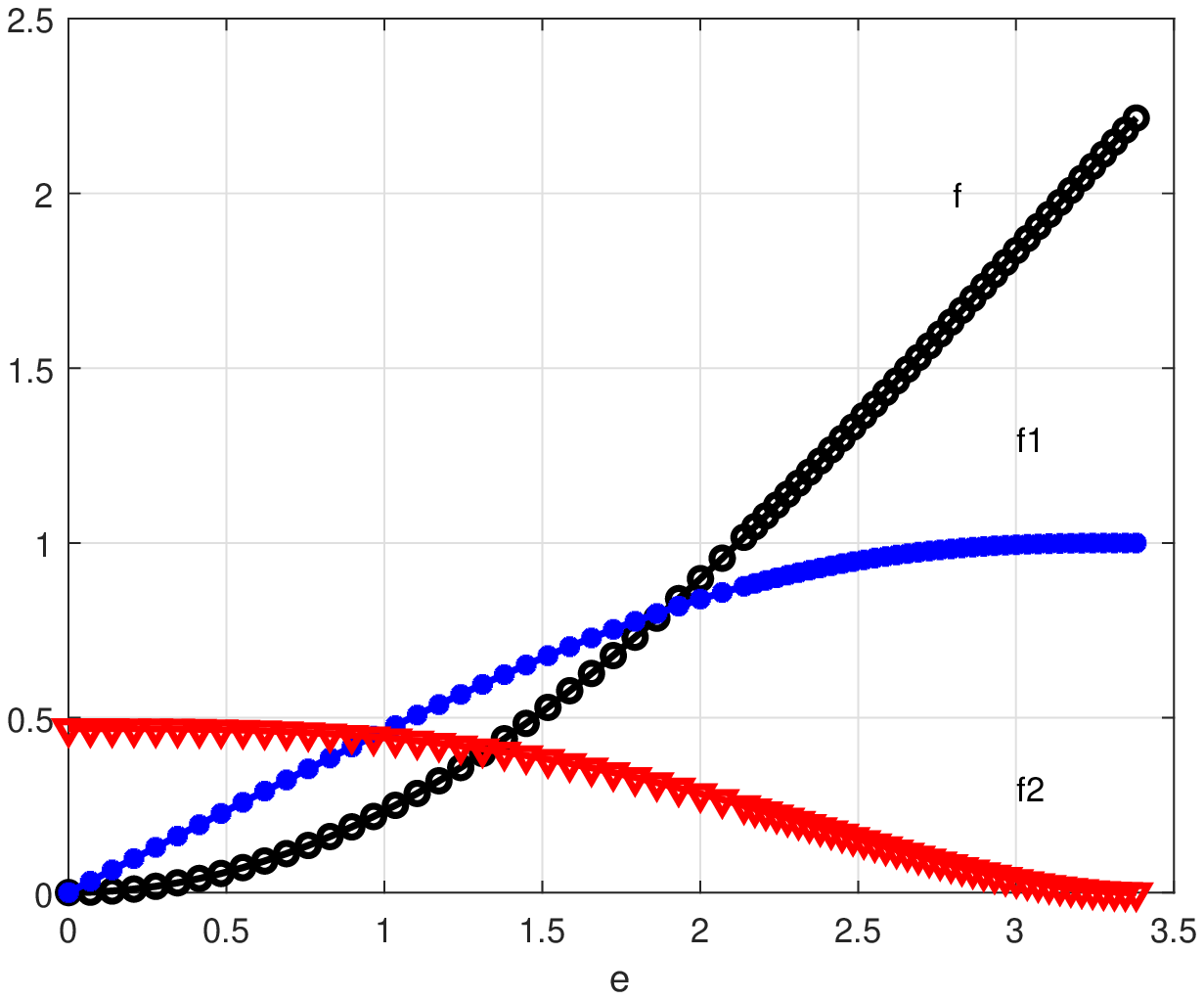}
\caption{Numerical result obtained using the \FBF \ for (\ref{eq:ExBlasius1}) with $P =3/2$ and $\e = 1{\D}-06$.}
	\label{fig:ExBlasius1}
\end{figure}

As far as the first extended Blasius problem is concerned, 
in table \ref{tab:NITM:missingIC} we list the chosen values of $\e$, the corresponding free boundary values $\eta_\e$, and the related missing initial conditions $\frac{d^2f}{d\eta^2}(0)$ for the problem (\ref{eq:ExBlasius1}) with $P = 3/2$.
\begin{table}[!hbt]
\caption{Numerical data and results.}
\vspace{.5cm}
\renewcommand\arraystretch{1.3}
	\centering
		\begin{tabular}{lr@{.}lr@{.}l}
\hline \\[-2.2ex]
{$\e$} & \multicolumn{2}{c}%
{$\eta_\e$} 
& \multicolumn{2}{c}%
{$ {\displaystyle \frac{d^2f}{d\eta^2}(0)}$}\\[1.5ex]
\hline
0.1      & 2 & 708708 & 0 & 482527634 \\ 
0.01     & 3 & 193357 & 0 & 469356138 \\
0.0001   & 3 & 323660 & 0 & 469098357 \\
0.00001  & 3 & 364091 & 0 & 469055438 \\
0.000001 & 3 & 376487 & 0 & 469055050 \\
0.0000001  & 3 & 380402 & 0 & 469055086 \\
0.00000001  & 3 & 381636 & 0 & 469055082 \\
0.000000001  & 3 & 382027 & 0 & 469055080 \\
0.0000000001 & 3 & 382150 & 0 & 469055080 \\
\hline			
		\end{tabular}
	\label{tab:NITM:missingIC}
\end{table}
The obtained value of the missing initial condition can be compared with the one $0.46905520505$ obtained by Fazio, using an iterative transformation method and reported in a recent preprint \cite{Fazio:2020:ITM}. 
As it is easily seen, the two values agree up to the first six decimal places.
Figure \ref{fig:ExBlasius1} shows the solution of the extended Blasius problem (\ref{eq:ExBlasius1}) with $P = 3/2$

\subsection{Second extended Blasius}
Once again, we rewrite the \FBF \ of the second extended Blasius problem in standard form.
We define the same variables as before: $ \theta = \eta/\eta_\e$, $u_1 = f(\eta)$, $u_2 = \frac{df}{d\eta}(\eta)$,$u_3 = \frac{d^2f}{d\eta^2}(\eta)$ and $u_4 =\eta_\e$,
so that, the \FBF \ for the problem (\ref{eq:ExBlasius2}) is
\begin{align}\label{eq:ExB2FBFNF} 
& {\displaystyle \frac{du_1}{d\theta}} = \eta_\e u_2 \nonumber \\
& {\displaystyle \frac{du_2}{d\theta}} = \eta_\e u_3 \nonumber \\
& {\displaystyle \frac{du_3}{d\theta}} = - \eta_\e \frac{u_1 u_3}{(P+1)*((P-1)*\left|u_3\right|^{P-2}+\left|u_3\right|^{P-1})}   \nonumber \\[-1.5ex]
& \\[-1.5ex]
& {\displaystyle \frac{du_4}{d\theta}} = 0 \nonumber\\
& f(0) = {\displaystyle \frac{df}{d\theta}}(0) = 0, \qquad
{\displaystyle \frac{df}{d\theta}}(1) = 1 \quad {\displaystyle \frac{d^2f}{d\theta^2}}(1) = \e \ , \nonumber 
\end{align}
For the numerical solution of the BVP (\ref{eq:ExB2FBFNF}) we used the \textit{bvp4c.m} MATLAB routine with initial iterate again provided by the relations (\ref{eq:IT}).

\begin{figure}[p]
\centering
\psfrag{e}[1][]{$\eta$}  
\psfrag{f1}[l][]{$\ds \frac{df}{d\eta}$} 
\psfrag{f2}[l][]{$\ds \frac{d^2f}{d\eta^2}$} 
\includegraphics[width=12cm,height=9cm]{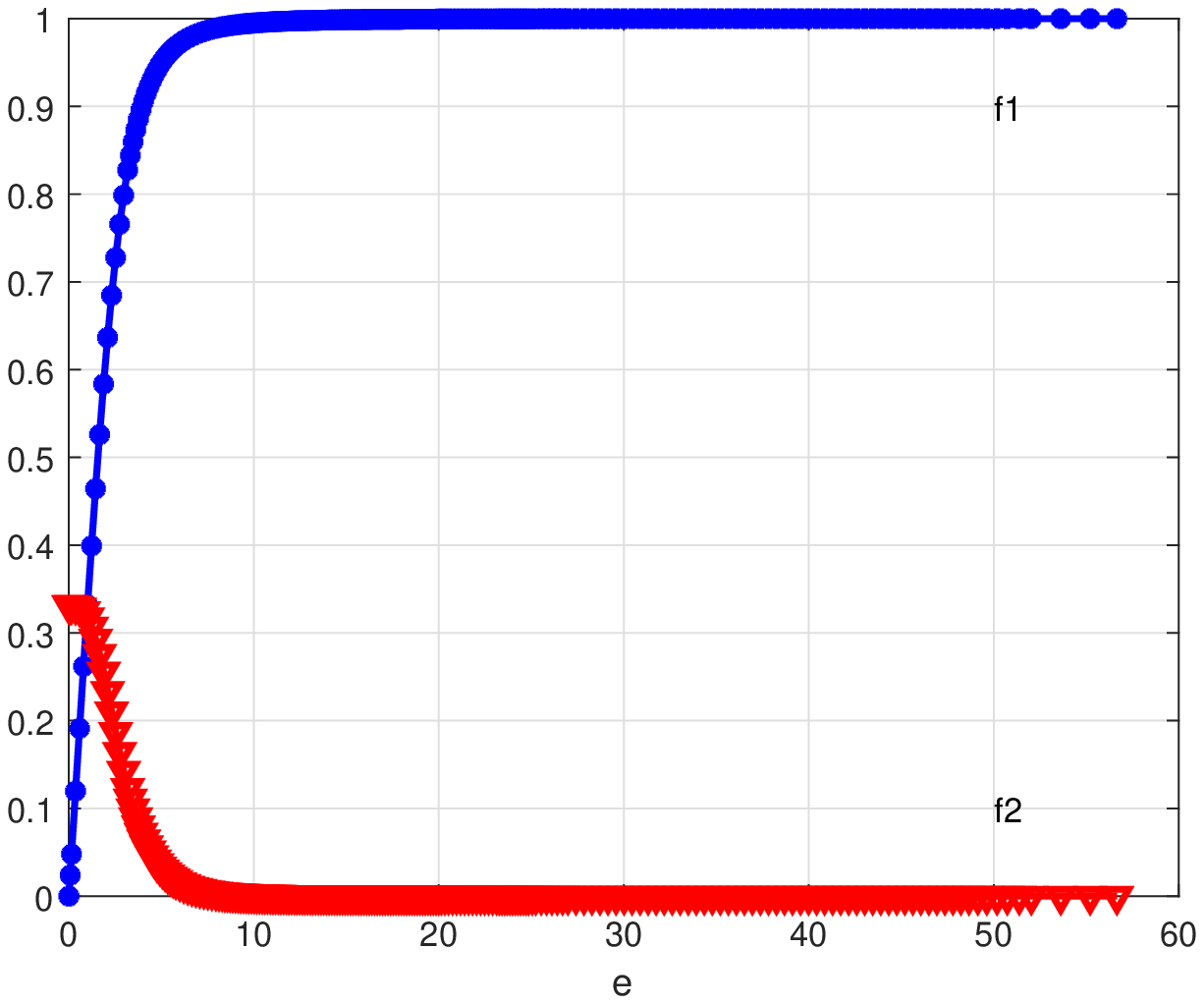} \\[1ex]
\includegraphics[width=12cm,height=9cm]{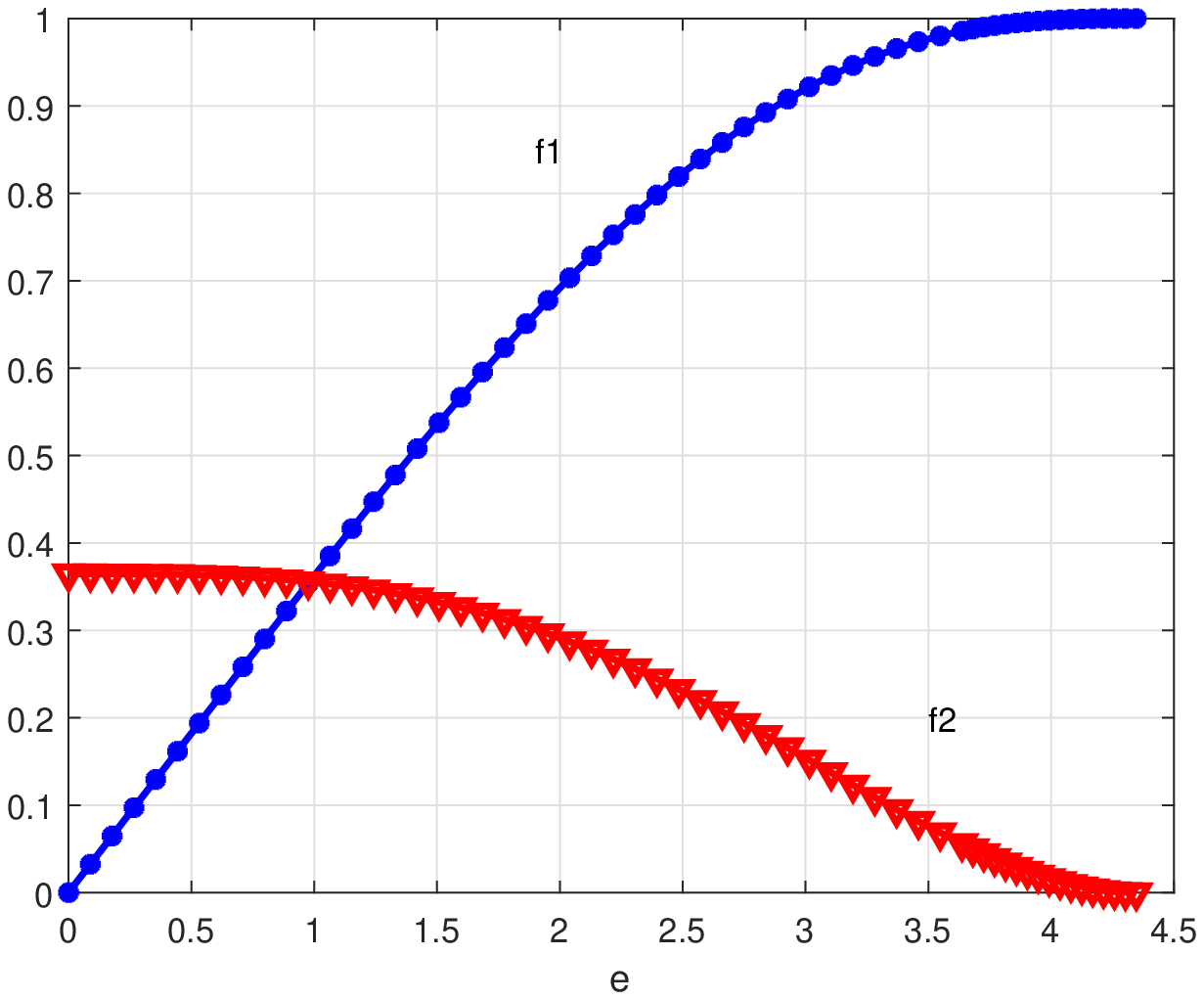}
\caption{Numerical result obtained using the \FBF \ with $\e = 1{\D}-06$ for (\ref{eq:ExBlasius2}) with: top frame $P =1/2$ and bottom $P = 3/2$.}
	\label{fig:ExBlasius2}
\end{figure}

Here, we report the numerical results obtained for the \FBF \ of the second extension (\ref{eq:ExBlasius2}) of the Blasius problem for two values of the involved parameter. 
The first value of the parameter is $P =1/2$, and in this case for $\e = 1{\D}-06$ we find $\eta_\e = 56.654480$ and $ \frac{d^2f}{d\eta^2}(0) = 0.331237479$.
The second considered value of the parameter is $P = 2$, and in this case for $\e = 1{\D}-06$ we find $\eta_\e = 4.346478$ and $ \frac{d^2f}{d\eta^2}(0) = 0.364773537$.
These values of the missing initial condition can be compared with those found by Fazio \cite{Fazio:2020:NIT}, using a non-iterative transformation method.
Those values are $0.337170$ for $P =1/2$ and $0.364772$ for $P=3/2$.
In figure \ref{fig:ExBlasius2}, for the reader convenience, we plot the two numerical solution components $f(\eta)$ and $\frac{df}{d\eta}(\eta)$.

\section{Conclusions} 
In this paper we have defined the \FBF \ for two extended Blasius problems.
These problems are of interest in boundary layer theory and are deduced from the governing partial differential equations by using appropriate similarity variables.
As far as the scaling invariance theory is concerned, we refer the interested reader to the books by Bluman and Cole \cite{Bluman:1974:SMD}, Barenblatt \cite{Barenblatt:1979:SSS}, or Dresner \cite{Dresner:1983:SSN}.

For both the extended Blasius problems, the computed results, for the so-called missing initial condition, are favourably compared with recent results available in the literature.

\vspace{1.5cm}

\noindent {\bf Acknowledgement.} {The research of this work was 
partially supported by the FFABR grant of the University of Messina and by the GNCS of INDAM.}


\begin{thebibliography}{10}

\bibitem{Akcay:1999:DRN}
M.~Akcay and M.~A. Y\"ukselen.
\newblock Drag reduction of a non-{N}ewtonian fluid by fluid injection on a
  moving wall.
\newblock {\em Arch. Appl. Mech.}, 69:215--225, 1999.

\bibitem{Ascher:1981:RBV}
U.~M. Ascher and R.~D. Russell.
\newblock Reformulation of boundary value problems into \lq \lq standard\rq \rq
  \ form.
\newblock {\em SIAM Rev.}, 23:238--254, 1981.

\bibitem{Barenblatt:1979:SSS}
G.~I. Barenblatt.
\newblock {\em Similarity, Self-Similarity and Intermediate Asymptotics}.
\newblock Consultant Bureau, New York, 1979.

\bibitem{Benlahsen:2008:GBE}
M.~Benlahsen, M.~Guedda, and R.~Kersner.
\newblock The generalized {B}lasius equation revisited.
\newblock {\em Math. Computer Model.}, 47:1063--1076, 2008.

\bibitem{Blasius:1908:GFK}
H.~Blasius.
\newblock Grenzschichten in {F}l\"{u}ssigkeiten mit kleiner {R}eibung.
\newblock {\em Z. Math. Phys.}, 56:1--37, 1908.

\bibitem{Bluman:1974:SMD}
G.~W. Bluman and J.~D. Cole.
\newblock {\em Similarity Methods for Differential Equations}.
\newblock Springer, Berlin, 1974.

\bibitem{Collatz}
L.~Collatz.
\newblock {\em The Numerical Treatment of Differential Equations}.
\newblock Springer, Berlin, 3rd edition, 1960.

\bibitem{deHoog:1980:ATB}
F.~R. de~Hoog and R.~Weiss.
\newblock An approximation theory for boundary value problems on infinite
  intervals.
\newblock {\em Computing}, 24:227--239, 1980.

\bibitem{Dresner:1983:SSN}
L.~Dresner.
\newblock {\em Similarity Solutions of Non-linear Partial Differential
  Equations}, volume~88 of {\em Research Notes in Math.}
\newblock Pitman, London, 1983.

\bibitem{Fazio:1992:BPF}
R.~Fazio.
\newblock The {Blasius} problem formulated as a free boundary value problem.
\newblock {\em Acta Mech.}, 95:1--7, 1992.

\bibitem{Fazio:1996:NAN}
R.~Fazio.
\newblock A novel approach to the numerical solution of boundary value problems
  on infinite intervals.
\newblock {\em SIAM J. Numer. Anal.}, 33:1473--1483, 1996.

\bibitem{Fazio:2020:NIT}
R.~Fazio.
\newblock A non-iterative transformation method for an extended {B}lasius
  problem, 2020.
\newblock Preprint available at the URL:
  \url{http://mat521.unime.it/~fazio/preprints/ExBlasius2020.pdf}.

\bibitem{Fazio:2020:ITM}
R.~Fazio and A.~Insana.
\newblock An iterative transformation method for a similarity boundary layer
  model, 2020.
\newblock Preprint available at the URL:
  \url{http://mat521.unime.it/~fazio/preprints/Fazio-Insana-2020R2.pdf}.

\bibitem{Fox}
L.~Fox.
\newblock {\em Numerical Solution of Two-point Boundary Value Problems in
  Ordinary Differential Equations}.
\newblock Clarendon Press, Oxford, 1957.

\bibitem{Kim:1983:MHT}
H.~W. Kim, D.~R. Jeng, and K.~J. De{W}itt.
\newblock Momentum and heat transfer in power-law fluid flow over
  two-dimensional or axisymmetic bodies.
\newblock {\em Inter. J. Heat Mass Transfer}, 26:245--259, 1983.

\bibitem{Lee:1966:SSN}
S.~Y. Lee and W.~F. Ames.
\newblock Similarity solutions for non-{N}ewtonian fluids.
\newblock {\em A.I.Ch.E. J.}, 12:700--708, 1966.

\bibitem{Lentini:BVP:1980}
M.~Lentini and H.~B. Keller.
\newblock Boundary value problems on semi-infinite intervals and their
  numerical solutions.
\newblock {\em SIAM J. Numer. Anal.}, 17:577--604, 1980.

\bibitem{Liao:2005:CNP}
S.-J. Liao.
\newblock A challenging nonlinear problem for numerical techniques.
\newblock {\em J. Comput. Appl. Math.}, 181:467--472, 1997.

\bibitem{Lin:1979:LBL}
F.~N. Lin and S.~Y. Chern.
\newblock Laminar boundary-layer flow of non-{N}ewtonian fluid.
\newblock {\em Inter. J. Heat Mass Transfer}, 22:1323--1329, 1979.

\bibitem{Markowich:TAS:1982}
P.~A. Markowich.
\newblock A theory for the approximation of solution of boundary value problems
  on infinite intervals.
\newblock {\em SIAM J. Math. Anal.}, 13:484--513, 1982.

\bibitem{Markowich:ABV:1983}
P.~A. Markowich.
\newblock Analysis of boundary value problems on infinite intervals.
\newblock {\em SIAM J. Math. Anal.}, 14:11--37, 1983.

\bibitem{Rubel:1955:EET}
L.~A. Rubel.
\newblock An estimation of the error due to the truncated boundary in the
  numerical solution of the {B}lasius equation.
\newblock {\em Quart. Appl. Math.}, 13:203--206, 1955.

\bibitem{Schlichting:2000:BLT}
H.~Schlichting and K.~Gersten.
\newblock {\em Boundary Layer Theory}.
\newblock Springer, Berlin, 8th edition, 2000.

\bibitem{Schowalter:1960:ABL}
W.~R. Schowalter.
\newblock The application of boundary-layer theory to power-law fluids: similar
  solutions.
\newblock {\em A.I.Ch.E. j.}, 6:24--28, 1960.

\end{thebibliography}
\end{document}